\documentclass{amsart}
\usepackage{amsmath,amsthm,amssymb,IMjournal}
\begin{document}
\newtheorem{theoreme}{Theorem}[section]
\newtheorem{lemme}[theoreme]{Lemma}
\newtheorem{proposition}[theoreme]{Proposition}
\newtheorem{corollaire}[theoreme]{Corollary}
\newtheorem{definition}[theoreme]{D\'efinition\rm}
\newtheorem{remarque}{\bf Remark}
\newtheorem{exemple}{\it Example\/}
\renewcommand\theequation{\arabic{equation}}

\newcommand{\ds}{\displaystyle}
\newcommand{\DD}{\mathbb D}
\newcommand{\TT}{\mathbb T}
\newcommand{\RR}{\mathbb R}
\newcommand{\NN}{\mathbb N}
\newcommand{\CC}{\mathbb C}

\newtheorem{The}{Theorem}[section]

\numberwithin{equation}{section}

\title{Estimates in the Hardy-Sobolev space of the  annulus and stability result}

\author{\|Imed |Feki|, Department of Mathematics, Faculty of Sciences, Sfax University,
 B.P 1171, Sfax 3018, Tunisia}



\abstract
The main purpose of this work is to establish some logarithmic estimates of optimal type in the  Hardy-Sobolev space $H^{k, \infty}; k \in {\mathbb{N}}^*$ of an annular domain. These results are considered as a continuation of a previous study in the setting of the unit disk by  L. Baratchart and M. Zerner: On the recovery of functions from pointwise boundary values in a Hardy-sobolev class of the disk. J.Comput.Apll.Math 46(1993), 255-69 and by S. Chaabane and I. Feki: Logarithmic stability estimates in Hardy-Sobolev spaces $H^{k,\infty}$. C.R. Acad. Sci. Paris, Ser. I 347(2009), 1001-1006.

As an application, we prove a logarithmic stability result for the inverse problem  of identifying a Robin  parameter  on a part of the boundary of an annular domain starting from its behavior on the complementary boundary part.
\endabstract



\thanks
This research has been supported by the Laboratory of Applied Mathematics and Harmonic Analysis: L. A. M. H. A. LR 11ES52
\endthanks

\section{Introduction}
The purpose of this paper is to establish  logarithmic  estimates of optimal type in the Hardy-Sobolev space $H^{1, \infty}(G_s)$ where $ s \in ]0, 1[$ and $G_s$ is the annulus of radius $(s,1)$. More precisely we study the behavior on the whole boundary of the annulus $G_s $ with respect to the uniform norm of any function $f$ in the unit ball of the Hardy-Sobolev space $H^{1, \infty}(G_s)$  starting from its behavior on any open  connected subset $I \subset \partial G_s$ with respect to the $L^1$-norm. Our result can be viewed as an extension of those established in \cite{chaabane fellah jaoua leblond 2004, Leblond, meftahi}.

The particular case where $I=\TT$ has been  considered   by L. Leblond, M. Mahjoub and J. R. Partington in \cite{Leblond}. The authors proved in this case that the $L^2$-norm of any function $f$ in the unit ball of the Hardy-Sobolev space $H^{1,2}(G_s)$ on the inner boundary $s\TT$  is controlled by the corresponding norm taken on the outer boundary $\TT$. In the same context, H. Meftahi and F. Wieolonsky gave recently in \cite{meftahi}  an explicit logarithmic inequality  exhibiting  the dependence with respect to the inner radius $s$ of the above control.
The first estimate of this kind remounts to L. Baratchart and M. Zerner where the authors proved in \cite{Baratchart}, a $\ds log log/log$ control with $L^2$-norm in the Hardy-Sobolev space $H^{1, 2}$ of the unit disk $\DD$. In \cite{Alessandrini piero rondi}, Alessandrini and al. have proved with quite different method an estimate of $\ds 1/log^\alpha$-type, $0< \alpha < 1$. Recently, the author of this paper together with S. Chaabane  \cite{Chaabane feki 2009} proved in the uniform norm some optimal logarithmic estimates in the Hardy-Sobolev space $H^{k, \infty}({\DD}); \ k \in \mathbb{N}^*$.

For more regular functions, we improve  inequality (\ref{eqn:main1}) amid the class of bounded $H^{k, \infty}(G_s)$ functions.

These logarithmic estimates allow us to prove stability result for the inverse problem of recovering a Robin coefficient on a part of the boundary of an annular domain starting from its behavior on the complementary boundary part. The particular case where the inaccessible part of the boundary is the  inner circle  has been proved in \cite{Leblond}. We can also refer the reader to \cite{Alessandrini piero rondi, Chaabane feki 2009, chaabane fellah jaoua leblond 2004, Leblond, meftahi} for stability estimates in the case of simply or doubly connected  domains.

\section {Notation and preliminary results}
Let $\DD$ be the open unit disk in $\mathbb{C}$ with boundary
$\TT$ and let  $G_s$ denote the annulus:
\begin{eqnarray}
G_s= \{ z \in {\mathbb{C}}; \quad s < |z| < 1\};\qquad 0< s < 1.\nonumber
\end{eqnarray}
The boundary of  the annular domain $G_s$ consists of two pieces $ s\TT$ and $\TT$:
$\partial G_s = s \TT \cup \TT$. Let $I$ be
 any connected open subset of the boundary of $G_s $ and let $ J = \partial G_s \setminus I$.
 We also equip the boundary $ \partial G_s $ with the usual  Lebesgue measure $\mu$ normalized so that the
 circles $\TT$ and $s\TT$, each have unit measure. Furthermore, we denote
 by $\ds\lambda = \mu(I)/(2\pi)$,  we  assume  that $\lambda \in  \ ]0, 1[$ and we defined by:
 \begin{eqnarray}
||f||_{L^1(I)} = \frac{1}{2\pi\lambda} \int_I |f(r e^{i\theta})| d \theta ,\nonumber
\end{eqnarray}
for the $L^1$-norm of $f$ on $I$, where $r=s$ if $I \subset s\TT$ and $r=1$ if $I \subset \TT$.

In the sequel, the Hardy space $ H^{\infty}(G_s)$  is defined
as the space of bounded analytic functions on $G_s$. According to (\cite{chevreau}, Theorem 7.1),
the Hardy space $ H^{\infty}(G_s)$ can be identified to the direct sum:
\begin{eqnarray}
H^{\infty}(G_s) = H^\infty({\mathbb{D}}) \oplus H^\infty_0({\mathbb{C}} \setminus s \overline {\mathbb{D}}),\nonumber
\end{eqnarray}
where the Hardy space $H^\infty_0({\mathbb{C}} \setminus s \overline {\mathbb{D}})$ is defined as the set of analytic functions in ${\mathbb{C}} \setminus s \overline {\mathbb{D}}$, with a zero limit at infinity.
Hence we can regard it as a closed subspace $ H^{\infty}(\partial G_s)$ of $ L^{\infty}(\partial G_s)$.
Equivalent definitions of Hardy spaces on annular domains are discussed by several authors
(\cite{BaratLeblondPartington, chalendar, chevreau, Rudin, sarason}). We can also refer the reader to \cite{Duren} for a more comprehensive details on Hardy spaces.

For $ k \in {\mathbb{N}}^*$, we  designate by $H^{k,\infty}(G_s)$, the Hardy-Sobolev space of order $k$
of the annulus:
\begin{eqnarray}
 \ds H^{k,\infty}(G_s)= \{f \in H^\infty(G_s) \ /\quad f^{(j)} \in
H^\infty(G_s),\quad j=0,...,k\},\nonumber
\end{eqnarray}
where $ f^{(j)}$ denotes the  $ j^{th}$ complex derivative of $ f $.
We endow $ H^{k,\infty}(G_s)$ with the norm inherited from the space $L^\infty(\partial G_s)$:
\begin{eqnarray}
 \ds ||f||_{{H}^{k,\infty}(G_s)} = {\ds \max_{0 \leq j \leq
k}\left(||f^{(j)}||_{L^\infty(s\TT)} \ +\ ||f^{(j)}||_{L^\infty(\TT)}\right)}. \nonumber
\end{eqnarray}
Let ${\mathcal{B}}_{k,\infty} = \left\{ f \in H^{k,\infty}(G_s); \;
\;||f||_{H^{k,\infty}(G_s)}\leq 1 \right\}$ be the closed unit ball of $ H^{k,\infty}$.

Next, we  introduce the Poisson kernel $p$  for the annulus $G_s$. Following  Sarason \cite{sarason} and Hwai \cite{hwai}, we consider the following holomorphic function:
\begin{eqnarray}
 \ds F(t,r ) = \frac{1}{2q_0}\ tanh\left(-\frac{\pi t}{2q_0} +\ i(\frac{\pi}{4} +\frac{\pi}{2q_0}log\frac{r}{\sqrt{s}})\right), \nonumber
\end{eqnarray}
where \ $q_0 = -log s,\quad 0< s < r <1 \quad\hbox{ and}\quad  t \in {\mathbb{R}}$.

The imaginary part $P(t,r)$ of $F(t,r )$ is the harmonic function given by:
\begin{eqnarray}
 \ds P(t,r) = \frac{1}{2q_0}\ \frac{ cos\left(\frac{\pi}{q_0}\ log \frac{r}{\sqrt{s}}\right)}
{cosh \frac{\pi t}{ q_0}\ -\ \sin\left(\frac{\pi }{ q_0}\ log \frac{r}{\sqrt{s}}\right)}.\nonumber
\end{eqnarray}
Referring to \cite[p. 92]{hwai}, we recall the following lemma.
\begin{lemme}\label{Lemma:kernelpropriety1}  The harmonic function $P$ satisfies the following properties:
\begin{itemize}
\item [i) ]
$\ds P(t,r) > 0 \quad \hbox{for} \quad s < r < 1\quad\hbox{ and}\quad  t \in {\mathbb{R}}$.\\
\item [ii) ]
$\ds \int_{-\infty}^{+\infty} P(t,r)\ dt \ +\  \int_{-\infty}^{+\infty} P(t,\frac{s}{r})\ dt = 1
\quad \hbox{for} \quad s < r < 1$.\\
\item [iii) ] There exists a non negative constant $C$\ such that for every $|t| \leq \pi$ and  $j$ large enough, we have:
\begin{eqnarray}
 \ds |P(t + 2\pi j, r)| \ \leq \ min \left(\frac{C}{j^4},\ \frac{C\ cos(\frac{\pi}{q_0})\ log(\frac{r}{\sqrt{s}})}{t^4}\right).\nonumber
\end{eqnarray}
\end{itemize}
\end{lemme}
This lemma allows us to define the Poisson kernel $p$ for the annular domain $G_s$:
\begin{eqnarray}
 \ds p(t,r) = \sum_{j =-\infty}^{+\infty} \ P(t + 2\pi j, r)\quad \hbox{for}\quad  |t| \leq \pi\quad \hbox{and}\quad  s < r < 1.\nonumber
\end{eqnarray}
We also have  from \cite{hwai},  the following lemma.
\begin{lemme}\label{Lemma:kernelpropriety2}\hskip10cm
\begin{itemize}
\item [i) ]
$\ds p(t,r)$   is a harmonic function on the annulus $ G_s$.\\
\item [ii) ]
$\ds p(t,r) > 0 \quad  \hbox{ for } \quad s < r < 1 \quad \hbox{and}\quad   |t| \leq \pi$.\\
\item [iii) ]
$\ds \frac{1}{2\pi}\int_{0}^{2\pi} p(t,r)\ dt \ +\  \frac{1}{2\pi}\int_{0}^{2\pi} p(t,\frac{s}{r})\ dt = 1
\quad \hbox{for} \quad  s< r < 1$.
\end{itemize}
\end{lemme}
In the next lemma, we recall the Poisson-Jensen formula for the annulus, see (\cite[p.25]{sarason}). This  will be of interest later.
\begin{lemme}\label{Lemma:logf}
Let $ f \not\equiv 0$ be a function in $H^q(G_s)$ for $ 1\leq q \leq \infty$. Then for all $r e^{it} \in G_s$, we have
\begin{eqnarray}
 \ds log|f(r e^{it})|\ \leq \ \frac{1}{2\pi}\int_{0}^{2\pi} p(t,r)log|f( e^{it})|\ dt +
\frac{1}{2\pi}\int_{0}^{2\pi} p(t,\frac{s}{r})log|f(s e^{it})|\ dt.\nonumber
\end{eqnarray}
\end{lemme}
\section{Optimal logarithmic estimates in  $\bf  H^{k,\infty}$ }
Our objective in this section, is to establish some logarithmic estimates in the Hardy-Sobolev space $H^{k, \infty}(G_s);\ k \in \mathbb{N}^*$ that can be viewed as a continuation of  the results  already established by  \cite{Chaabane feki 2009,  Leblond, meftahi}. We start by recording a variant of the Hardy-Landau-Littlewood inequality which will crucially be used in the proofs of Theorem \ref{theorem:fundamentel} and Theorem \ref{theorem:kinfty}, see \cite[chapter VIII p.147]{Brezis} and \cite{Nirenberg}.
\begin{lemme}\label{lemma:HLL}
Let $\mathcal{I}$ be a bounded interval and let $j \in {\NN}$\ such that $j \geq 2 $. Then,
there exists a non negative constants $C_\infty(\mathcal{I}, j)$\ such that
\begin{eqnarray}\label{HLL}
 {\ds
||g'||_{L^\infty({\mathcal{I}})} \leq C_\infty({\mathcal{I}}, j)
||g||^{1/j}_{W^{j, \infty}({\mathcal{I}})} ||g||^{1-1/j}_{L^\infty({\mathcal{I}})}} \quad \hbox{for all}\  g \ \in W^{j, \infty}({\mathcal{I}}).
\end{eqnarray}
\end{lemme}
Next, we give a lower bound for the  Poisson kernel $p$   which will be useful for the proof of Lemma \ref{lemma:BZ}.
\begin{lemme}\label{Lemma:pminoration}
There exists a non negative constant $C_s$ depending only on $s \in ]0, 1[$, such that for every $|t|\leq \pi$, we have
\begin{eqnarray}
 \begin{array}{cclcr}
 \ds  p(t,r) &\ \geq &\frac{2 C_s}{log\ s}\ (log\ s \  -\ log\ r)& \hbox {if }& s < r \leq \sqrt{s}.\\
&&&&\\
 \ds p(t,r) &\ \geq &\frac{2 C_s}{log\ s}\ log\ r &\hbox {if }& \sqrt{s} \leq r < 1.
\end{array}\nonumber
\end{eqnarray}
\end{lemme}
{\bf{Proof.}} Remind the reader that $q_0 = -log s$.  Let $r \in ]s, 1[$, then \  $\ds \frac{\pi}{q_0}\ log (\frac{r}{\sqrt{s}}) \ \in \ ]-\frac{\pi}{2}, \ \frac{\pi}{2}[$  and therefore
\begin{eqnarray}
 \ds P(t+2\pi j, r) \ \geq \ \frac{1}{2q_0\left(1 + cosh \frac{\pi}{q_0}(t + 2 \pi j)\right)}
\ cos\left( \frac{\pi}{q_0}\ log (\frac{r}{\sqrt{s}})\right).\nonumber
\end{eqnarray}
Since
\begin{eqnarray}
 \ds C_s(t) = \frac{1}{2q_0}\ \sum_{j =-\infty}^{+\infty} \ \frac{1}{ 1 + cosh \frac{\pi}{q_0}(t + 2 \pi j)} \ < \ \infty \qquad \hbox {for every } \ |t|\leq \pi,\nonumber
\end{eqnarray}
we deduce that
\begin{eqnarray}\label{minorationdep}
  p(t, r) \ \geq \ C_s \ cos\left(\frac{\pi}{q_0}\ log (\frac{r}{\sqrt{s}})\right), \quad \ds C_s = \inf_{|t|\leq \pi} \ C_s(t).
\end{eqnarray}
In the case where  $ r \in ]s, \sqrt{s}]$, we have
$\frac{\pi}{q_0}\ log (\frac{r}{\sqrt{s}}) \in \ ]-\frac{\pi}{2}, 0]$.
Using the inequality \   $ cos \ x \geq \frac{2}{\pi}\  x +1 $ \ for  \  $x\  \in \ ]-\frac{\pi}{2},  0]$, we obtain
\begin{eqnarray}
  p(t,r) \ \geq \frac{2 C_s}{log\ s}\ (log\ s \  -\ log\ r). \nonumber
 \end{eqnarray}
Otherwise,  $ r \in [\sqrt{s}, 1[$ and
$\frac{\pi}{q_0}\ log (\frac{r}{\sqrt{s}}) \in [0,  \frac{\pi}{2}[.$
Using the inequality \ $ cos \ x \geq -\frac{2}{\pi} x +1$ \  for  \  $x \ \in \ [0,  \frac{\pi}{2}[$, we obtain
\begin{eqnarray}
 \ds p(t,r) \ \geq \frac{2 C_s}{log\ s}\  log\ r,
 \end{eqnarray}
 which achieves the proof of the Lemma.

We adapt the same arguments developed in (\cite{Baratchart}, lemma 4.1)
with some slight shifts to prove the following
 \begin{lemme}\label{lemma:BZ}
Let $ g \in  H^{\infty}(G_s)$ and $\displaystyle m \geq  ||g||_{L^\infty(\partial G_s)}$. Then, for every
$z \in\overline{G}_s$, we have
\begin{eqnarray}
 \begin{array}{lllcclcr}
 |g(z)| & \leq &
 \ds  m \ \left\|\frac{g}{m}\right\|_{L^1(I)}^{\frac{2 \lambda C_s}{log s}\ (log s\ -\ log|z|)}&
\hbox {\rm if }&  s & < |z| &\leq &\sqrt{s},\\
&&&&&\\
 \ds  |g(z)| & \leq &  m\  \left\|\frac{g}{m}\right\|_{L^1(I)}^{\frac{2 \lambda C_s }{log s}\ log|z|}&
\hbox{ \rm if }& \sqrt{s}& \leq |z|& < &1.\\
  \end{array}\nonumber
\end{eqnarray}
 \end{lemme}
{\bf{Proof.}}
Let $ h = g/m$ and let $z = r e^{it} \in G_s$. From Lemma \ref{Lemma:logf} and the fact that $log|h|$ is a non positive subharmonic function, we get
\begin{eqnarray}
 log(|h(r e^{it})|)  \leq  \frac{1}{2\pi} \int_{0}^{2\pi}p(t-\theta , r)\ log(|h(e^{i\theta})|)\ d\theta
+ \frac{1}{2\pi} \int_{0}^{2\pi}p(t-\theta , \frac{s}{r})\ log(|h(se^{i\theta})|)\ d\theta.\nonumber
\end{eqnarray}
If we suppose  that $I \subset \TT$, then by using the facts that $p(t,r) > 0$ and that $log|h|\ \leq \ 0$,  we deduce that
\begin{eqnarray}
 log(|h(r e^{it})|)  \leq  \lambda \int_I p(t-\theta , r)\ log(|h(e^{i\theta})|)\ \frac{d\theta}{2\pi\lambda},\nonumber
\end{eqnarray}
consequently, from Lemma \ref{Lemma:pminoration}, we obtain
$$
\begin{array}{rrllll}
&log\left( |h(z)|\right ) & \ds  \leq &\frac{2\lambda C_s}{log\ s}\ (log s  -  log r)\
 {\ds \int_I log\left(|h(e^{i\theta})|\right) \frac{d\theta}{2\pi\lambda}}
&\hbox{if}& s < |z| \leq \sqrt{s},\\
& & & &\\
&log\left( |h(z)|\right ) &\ds  \leq & \frac{2\lambda C_s}{log s}  \ log r
 {\ds \int_I log\left(|h(e^{i\theta})|\right)\ \frac{d\theta}{2\pi\lambda}}
& \hbox{if}& \sqrt{s} < |z| <1
\end{array}
$$
By using Jensen's inequality,  we deduce that
$$
\begin{array}{rrllll}
 &|g(z)| & \leq & \ds m \left\|\frac{g}{m}\right\|_{L^1(I)}^{\frac{2 \lambda C_s}{log s}\ (log s\ -\ log|z|)}
&\quad \hbox{if} & s < |z| \leq \sqrt{s},\qquad \qquad \qquad\quad\\
 &&&&\\
 &|g(z)| & \leq & \ds m \left\|\frac{g}{m}\right\|_{L^1(I)}^{\frac{2 \lambda C_s}{log s}\ log|z|}
& \quad \hbox{if} & \sqrt{s} < |z| <1.\qquad \qquad\qquad\quad
\end{array}
$$
If we suppose  that $I \subset s\TT$, then by using again the facts that $p(t,r) > 0$ and that $log|h|\ \leq \ 0$,  we get
\begin{eqnarray}
 log(|h(r e^{it})|)  \leq  \lambda \int_I p(t-\theta , \frac{s}{r})\ log(|h(se^{i\theta})|)\ \frac{d\theta}{2\pi\lambda}\nonumber
\end{eqnarray}
and the proof can be completed in a similar way as in the first case.

Let $f \in  H^{\infty}(G_s)$ and let $t$ be a real number such that  $|t| \leq \pi$. We designate by $F_t$ the radial primitive of $f$ that vanishes at $s$ and defined by:
\begin{eqnarray}\label{defprimitif}
 F_t(r) =\ds  \int_s^r f(x e^{it})\ dx \quad \hbox{for all}\quad   r \ \in\ {\mathcal{I}} = ]s, 1[.
\end{eqnarray}
From Lemma \ref{lemma:BZ}, we obtain
\begin{lemme}\label{lemma:primitive}
Let $f\in  H^{\infty}(G_s)$ and  $\displaystyle m \geq ||f||_{L^\infty(\partial G_s)}$. We suppose that $f$ is not identically zero and that \quad  $ ||f||_{L^1(I)} <  e^{-\frac{q_0}{\lambda C_s}}$. Then for all $|t| \leq \pi$ and
$ r  \in \ ]s, 1[$ we get

\begin{eqnarray}\label{radialprimitive}
 \ds |F_t(r)|\ \leq \ \ds \frac{(2s+1)q_0 m}{|2 \lambda C_s \ log\|\frac{f}{m}\|_{L^1(I)}|}
\end{eqnarray}
\end{lemme}
{\bf{Proof.}}
Let $ |t| \leq \pi$ and let  $ r \in ]s, 1[$.  From (\ref{defprimitif}) and the monotonicity of the function $\eta(y) = \int_s^y|f(xe^{it})| dx$, we have
\begin{eqnarray}
 |F_t(r)|  \leq \int_s^{\sqrt{s}}\ |f(x e^{it})|\ dx\ + \ \int_{\sqrt{s}}^1r\ |f(x e^{it})|\ dx,\nonumber
\end{eqnarray}
then according to Lemma \ref{lemma:BZ} we get
$$
\begin{array}{rrllll}
 & |F_t(r)| & \leq & m {\ds \int_s^{\sqrt{s}}\ \left\|\frac{f}{m}\right\|_{L^1(I)}^{\frac{2 \lambda C_s}{log s}(log s\ -\ log x)}\ dx}
  & + & m {\ds \int_{\sqrt{s}}^1\ \left\|\frac{f}{m}\right\|_{L^1(I)}^{\frac{2 \lambda C_s}{log s}\ log x}\ dx}\\
 &&&&&\\
 & & \leq & \ds \frac{m s}{|1 + \frac{2 \lambda C_s}{q_0} \ log\left\|\frac{f}{m}\right\|_{L^1(I)}|}
  & + & \ds \frac{m }{|1 - \frac{2 \lambda C_s}{q_0} \ log\left\|\frac{f}{m}\right\|_{L^1(I)}|}.
  \end{array}
$$
From the assumption that $ ||f||_{L^1(I)} <  e^{-\frac{q_0}{\lambda C_s}}$,  we have
\begin{eqnarray}\label{eqfm}
 \ds \frac{1}{| 1 + \frac{2 \lambda C_s}{q_0} \ log\|\frac{f}{m}\|_{L^1(I)}|}\ \  \leq \ \ds \frac{2}{|\frac{2 \lambda C_s}{q_0} \ log\|\frac{f}{m}\|_{L^1(I)}|},\nonumber
 \end{eqnarray}
and therefore, we conclude the desired inequality
\begin{eqnarray}
\ds  |F_t(r)|  \leq \frac{(2s+1) q_0m}{|2 \lambda C_s \ log\|\frac{f}{m}\|_{L^1(I)}|}.\nonumber
 \end{eqnarray}

We are now in a position to establish the main control theorem in the
Hardy-Sobolev space $ H^{1,\infty}(G_s)$.
\begin{theoreme}\label{theorem:fundamentel}
Let $f \in {\mathcal{B}}_{1,\infty}$ and  $\displaystyle m \geq ||f||_{L^\infty(\partial G_s)}$. We suppose that $f$ is not identically zero and that \quad  $ ||f||_{L^1(I)} <  e^{-\frac{q_0}{\lambda C_s}}$. Then
 \begin{eqnarray}\label{eqn:main1}
\|f \|_{L^\infty(\partial G_s)} \leq
\frac{C^2_\infty({\mathcal{I}}, 2)/(1 - 1/2e)}{|\lambda_0 Log\|f\|_{L^1(I)}|},
\end{eqnarray}
where  $\ds \lambda_0= min \left(1, \frac{2 \lambda C_s}{(1+2s)\ q_0}\right)$.

Moreover, for $I = {\mathbb{T}}$, there exists a sequence of functions $f_n \in {\mathcal{B}}_{1,\infty}$  such that
\begin{eqnarray}\label{eqn:sequence}
 {\ds \lim_{n \rightarrow +\infty} \|f_n
\|_{L^\infty(\partial G_s)}\left|Log\|f_n\|_{L^1({\mathbb{T}})}\right| \geq s|log s|.}
\end{eqnarray}
\end{theoreme}
{\bf Proof.} Let for every $|t| \leq \pi $, $F_t$ be the radial primitive of $f$ defined by equation (\ref{defprimitif})  and let $\displaystyle m \geq max(||f||_{L^\infty(\partial G_s)}, 1)$.
According to Lemma \ref{lemma:primitive}, we have
\begin{eqnarray}\label{eqn:prim0}
 |F_t(r)|\ \leq \ \frac{m}{|\lambda_0 \ log\|\frac{f}{m}\|_{L^1(I)}|}, \qquad \hbox{where}\quad
\lambda_0 =min (1, \frac{2 \lambda C_s}{(1+2s)\ q_0}).
\end{eqnarray}
Since $f \in {\mathcal{B}}_{1,\infty}$, then according to the Hardy-Landau-Littlewood inequality (\ref{HLL}),
 there exists a non negative constant $C= C_\infty({\mathcal{I}}, 2)$ such that
\begin{eqnarray}
  \|f\|_{L^\infty(\partial G_s)}\leq C \|F\|^{1/2}_{L^\infty(G_s)},\nonumber
\end{eqnarray}
and consequently,
\begin{eqnarray}\label{eqn:etapeinit}
 \|f\|_{L^\infty(\partial G_s)}\leq m_1:=C \left( \frac{m}{|\lambda_0 \ log\|\frac{f}{m}\|_{L^1(I)}|}\right)^{1/2}.
\end{eqnarray}
Making use of (\ref{eqn:prim0})  and (\ref{eqn:etapeinit}) for the new estimate
$m_1$ of $\|f\|_{L^\infty(\partial G_s)}$, one obtains
\begin{eqnarray}\label{eqn:etape1}
 \|f\|_{L^\infty(\partial G_s)} \leq \ C \left(\frac{m_1}{|\lambda_0 \ log\|\frac{f}{m_1}\|_{L^1(I)}|}\right)^{1/2}.
\end{eqnarray}
Let $\eta(x) = x|log x|^{1/2} $ and $\alpha = 1 - \frac{1}{2e}$. Since
$ m \geq 1,\ \lambda_0 \leq 1$ and $ g(x) \leq x^\alpha $ in $ ]0,  1] $, we get
\begin{eqnarray}
 {\ds \left\|\frac{f}{m_1}\right\|_{L^1(I)} =
\frac{(m \lambda_0)^{1/2}}{C} \ \eta \left(\left\|\frac{f}{m}\right\|_{L^1(I)}\right)\ \leq\
\|f\|_{L^1(I)}^{\alpha}}.\nonumber
\end{eqnarray}
From (\ref{eqn:etape1}) and the monotonicity of the mapping
$\ds\varepsilon(x)=\frac{1}{|Log \, x|}$, we obtain
\begin{eqnarray}
 \|f\|_{L^\infty(\partial G_s)}\leq C^{1+1/2} \
\frac{m^{(1/2)^2}{(\frac{1}{\alpha})^{1/2}}}{\left|\lambda_0
Log\|f\|_{L^1(I)}\right |^{1/2(1+1/2)}}.\nonumber
\end{eqnarray}
Proceeding thus repeatedly, we obtain for every $k \in \mathbb{N^*}$,
\begin{eqnarray}
 \|f\|_{L^\infty(\partial G_s)}\leq C^{b_k} \
\frac{m^{(1/2)^{k+1}}{(\frac{1}{\alpha})^{c_k}}}{\left|\lambda_0
Log\|f\|_{L^1(I)}\right |^{a_k}},\nonumber
\end{eqnarray}
where $a_k$, $b_k$ and $c_k$ are three recurrent sequences
satisfying
\begin{eqnarray}
 a_1 = \frac{1}{2}(1+\frac{1}{2}),\ b_1 = 1+\frac{1}{2} ,\
c_1 = \frac{1}{2}, \  a_{k+1} = \frac{1+a_k}{2},\  b_{k+1} =
1+\frac{b_k}{2}, \  c_{k+1} = \frac{1+c_k}{2}.\nonumber
\end{eqnarray}
The proof of inequality (\ref{eqn:main1}) is completed by letting $k
\rightarrow +\infty $.

To prove equation (\ref{eqn:sequence}), we consider the sequence of
functions, $ u_n(z) =\frac{1}{z^n};  \ n \in {\NN}^*.$

Let $I = {\mathbb{T}}$ and let
$ f_n\ =\ {\ds \frac{u_n}{\|u_n \|_{H^{1,\infty}(G_s)}}}$ be
the $H^{1,\infty}(G_s)$ normalized function of $u_n$. Then,
\begin{eqnarray}
 \|f_n \|_{L^\infty(s\mathbb{T})} \  = \frac{\frac{1}{s^n}}{n(1+ \frac{1}{s^{n+1}})},\
 \|f_n\|_{L^\infty(\mathbb{T})} = \frac{1}{n(1+ \frac{1}{s^{n+1}})} \quad \hbox{\rm
and}\  \|f_n \|_{L^\infty(\partial G_s)} \  = \frac{1+\frac{1}{s^n}}{n(1 + \frac{1}{s^{n+1}})}.\nonumber
\end{eqnarray}
Let $A_n = {\ds\|f_n\|_{L^\infty(\partial G_s)}\left|\ log\|f_n\|_{L^\infty(\mathbb{T})}\right|},$
then we have
\begin{eqnarray}
 A_n = s \frac{1+s^n}{n(1+s^{n+1})}\  | log n \ +\ log (1+s^{n+1})\ - \ (n+1)log s|.\nonumber
\end{eqnarray}
Hence, $ {\ds \lim_{n \rightarrow \infty}\ A_n = s |log s| }$ and this completes the proof.

\begin{remarque}
The estimate (\ref{eqn:main1}) still holds in more general situations of a smooth doubly-connected domain
$G \subset {\mathbb{R}}^2$(we can see \cite{Pommerenke} for more details on conformal mapping).
\end{remarque}
\begin{remarque}
The estimate (\ref{eqn:main1}) of Theorem \ref{theorem:fundamentel} is of optimal type: it is impossible to
find a function $\varepsilon$ which  tends to zero at zero such that for all $f \in {\mathcal{B}}_{1,\infty}$,
\begin{eqnarray}
  \|f \|_{L^\infty(\partial G_s)}\ \leq \
\frac{1}{|Log\|f\|_{L^1(I)}|}\varepsilon\left(||f||_{L^1(I)}\right).\nonumber
\end{eqnarray}
\end{remarque}
\begin{remarque}
The estimate (\ref{eqn:main1}) of Theorem \ref{theorem:fundamentel} is false in the general setting of bounded function $f \in H^\infty(G_s)$( we consider the $H^\infty$-normalized function of $u_n$).
\end{remarque}
\begin{remarque}
The question under investigations is to give  the optimal constant $C$ in equation  (\ref{eqn:main1}).
\begin{eqnarray}
 C = \ds \max_{f \in {\mathcal{B}}_{1,\infty}}\ \|f\|_{L^\infty(\partial G_s)}\left|Log\|f\|_{L^1(I)}\right|.\nonumber
\end{eqnarray}
\end{remarque}
The following corollary is a direct consequence of Theorem \ref{theorem:fundamentel}.
\begin{corollaire}
Let $K >0$ and  $f \in H^{1,\infty}(G_s)$ such that  $\|f
\|_{H^{1,\infty}(\partial G_s)} \leq K$ and \ $\|f \|_{L^1(I)} <  e^{-\frac{q_0}{\lambda\ C_s}}$. Then, we have
\begin{eqnarray}
  \|f \|_{L^\infty(\partial G_s)}\ \leq \
\frac{C^2_\infty({\mathcal{I}}, 2)\ max(1,K)/(1 - 1/2e)}{|\lambda_0 Log\|f\|_{L^1(I)}|}.\nonumber
\end{eqnarray}
\end{corollaire}
If we suppose that $f$ is more  regular, then we can improve
inequality (\ref{eqn:main1}) in the same way as in the proof of
Theorem \ref{theorem:fundamentel}.

\begin{theoreme}\label{theorem:kinfty}
Let $k \in \mathbb{N}^*$. There exists a non negative constant
$C$ depending only on $k, s$ and $\lambda$  such that for every $f \in {\mathcal{B}}_{k,\infty}$ satisfying
$ ||f||_{L^1(I)} <  e^{-\frac{q_0}{\lambda C_s}}$  also satisfies
\begin{eqnarray}
  \|f \|_{L^\infty(\partial G_s)}\ \leq \
\frac{C_k(s)}{| Log\|f\|_{L^1(I)}|^k}.\nonumber
\end{eqnarray}
Moreover, for $I = {\TT}$,  there exists a sequence $f_n$ of ${\mathcal{B}}_{k,\infty}$
such that
\begin{eqnarray}\label{eqn:sequencek}
 \ds \lim_{n \rightarrow +\infty} \|f_n
\|_{L^\infty(\mathbb{T})}\left|Log\|f_n\|_{L^1(I)}\right|^k \geq s|log s|^k.
\end{eqnarray}
\end{theoreme}

{\bf Proof.}  For every $|t| \leq \pi $, we have according to the proof of the previous theorem that  the radial primitive $F_t$ of $f$ satisfies inequality (\ref{eqn:prim0}).
Since $f \in {\mathcal{B}}_{k,\infty}$, then from  the Hardy-Landau-Littlewood inequality (\ref{HLL}) applied to $j = k+1$, we prove  that there exists a non negative constant $C= C_\infty({\mathcal{I}}, k+1)$ such that
\begin{equation}\label{eqn:etapeinitk}
 \|f\|_{L^\infty(\partial G_s)}\leq m_1:=C \left( \frac{m}{|\lambda_0 \ log\|\frac{f}{m}\|_{L^1(I)}|}\right)^{\frac{k}{k+1}}.
\end{equation}
Similarly to the proof of  Theorem \ref{theorem:fundamentel},  consider $ \rho = \frac{k}{k+1} ,\quad g_\rho(x) = x|log x|^{\rho} $ and $\sigma=1-\frac{\rho}{e}$. Then we have
 $ g_\rho(x) \leq x^\sigma $ in $ ]0 \, , \, 1]$ and consequently we establish
 for every $ j \in \mathbb{N}^*$ the following inequality
\begin{equation}
 \|f\|_{L^\infty(\partial G_s)}\leq C^{b_j} \
\frac{m^{(\rho)^{j+1}}{(\frac{1}{\sigma})^{c_j}}}{\left|\lambda_0
Log\|f\|_{L^1(I)}\right |^{a_j}},\nonumber
\end{equation}
where $a_j$, $b_j$ and $c_j$ are three recurrent sequences
satisfying
\begin{eqnarray}
 a_1 = \rho(1+\rho);\  b_1 = 1+\rho;\
c_1 = \rho;\  a_{j+1} = \rho (1+a_j);\  b_{j+1} =
1+\rho\ b_j \ \hbox{and} \  c_{j+1} = \rho (1+c_j).\nonumber
\end{eqnarray}
Then by letting $j \rightarrow +\infty $, we obtain
\begin{eqnarray}
 {\ds \|f\|_{L^\infty(\partial G_s)}\leq \frac{C_k(s)}{| Log\|f\|_{L^1(I)}|^k}}.\nonumber
\end{eqnarray}
To prove equation (\ref{eqn:sequencek}), we consider the same sequence as in the proof of  equation (\ref{eqn:sequence}), with the suitable $H^{k, \infty}(G_s)$ normalization norm.

\begin{corollaire}\label{cor:k-jinfty}
Let $K > 0, \ j$ and $k$\ be some integers with \ $0 \leq j < k$. Let  $f \in {H}^{k,\infty}$ such that $||f||_{H^{k, \infty}(G_s)} \leq K$ and $ ||f||_{H^{j,\infty}(I)} <  e^{-\frac{q_0}{\lambda C_s}}$.
  Then, there exist non negative constants
$C, \varepsilon$ depending only on $K, k, j, s$ and $\lambda$ such that
\begin{eqnarray}
  \|f \|_{H^{j, \infty}(\partial G_s)}\ \leq \
\frac{C}{| Log\|f\|_{L^1(I)}|^{k-j}}.\nonumber
\end{eqnarray}
provided that  $||f||_{L^1(I)} < \varepsilon$.
\end{corollaire}
{\bf Proof.}
Let $K_1 = \max(K, 1)$ and let $g = f/K_1$, then the derivative $g^{(i)}$ of order $i \in \{0,...,j\}$ belongs to  ${\mathcal{B}}_{k-i, \infty}$ and satisfy the assumptions of Theorem \ref{theorem:kinfty}. Hence, there exists a non-negative constant $C_1$ depending only on $K, k, i, s$ and $\lambda$ such that
\begin{eqnarray}\label{ineq-main-j}
    \|g^{(i)}\|_{L^\infty(\partial G_s)}\leq\frac{C_1}{|\log\|g^{(i)}\|_{L^1(I)}|^{k-i}}.
\end{eqnarray}
According to \cite[Theorem 1]{Nirenberg} and the assumption that $g \in {\mathcal{B}}_{k, \infty}$, there exist a non-negative constant $C_2$ such that
\begin{eqnarray}
\|g^{(i)}\|_{L^1(I)}\leq C_2\|g\|_{L^1(I)}^{1-i/k},\nonumber
\end{eqnarray}

we derive from (\ref{ineq-main-j}) and the monotonicity of the mapping $\eta_i(x) = 1/(log(1/x))^{k-i}$ that
\begin{eqnarray}\label{eqi-k}
\|g^{(i)}\|_{L^\infty(\partial G_s)}\leq C_1 \eta_i(C_2\|g\|_{L^1(I)}^{1-i/k}).
\end{eqnarray}
Let us choose $\varepsilon > 0$ small enough such that
\begin{eqnarray}\label{eqi-k2}
 \eta_i(C_2\|g\|_{L^1(I)}^{1-i/k})\leq 2 \eta_i(\|g\|_{L^1(I)}),
\end{eqnarray}
then from (\ref{eqi-k}) and (\ref{eqi-k2}), we obtain
\begin{eqnarray}
    \|g^{(i)}\|_{L^\infty(\partial G_s)}\leq\frac{2C_1}{|\log\|g\|_{L^1(I)}|^{k-i}}.\nonumber
\end{eqnarray}
Taking the maximum over all $i =0,...,j$ we achieved the proof of the corollary.

As an immediate consequence, we prove that if the $L^1$-norm of a bounded $H^{k, \infty}(\partial G_s)$ function is known to be small on a connected open  subset $I $ of $\partial G_s$ it remains also small (with uniform norm) on the whole boundary $\partial G_s$. The same result with $L^2$-norm has been established by Leblond and al. in \cite{Leblond}.
\begin{corollaire}
Let $j$ and $k$ be some integers with $0 \leq j < k$, and let $I \subset \partial G_s$ be
 any connected open subset. Let $(f_p)$  be a sequence of functions in the unit ball of the Hardy-Sobolev spaces $H^{k, \infty}(\partial G_s)$ such that $\|f_p\|_{L^1(I)} \longrightarrow 0$. Then $\|f_p\|_{H^{j,\infty}(\partial G_s)} \longrightarrow 0$.
\end{corollaire}

In the particular case where $I = {\TT}$, the following corollary provides logarithmic estimates with respect to the
 $L^\infty$-norm similar to those proved with the  $L^2$-norm by Leblond and al. in \cite{Leblond}.
\begin{corollaire}
Let $I = {\TT}$, $k$ and $j$ be some integers with  $0 \leq j < k$.  Then, there exist non negative constants
$C, \varepsilon$ depending only on $K, k, j$ and $I$ such that  whenever  $ f \in {\mathcal{B}}_{k, \infty}$  and satisfies    $ ||f||_{H^{j,\infty}(I)} <  e^{-\frac{q_0}{\lambda C_s}}$, we have
\begin{eqnarray}
  \|f \|_{H^{j,\infty}(s{\mathbb{T}})}\ \leq \
\frac{C} {| Log\|f\|_{L^1({\TT})}|^{k-j}}\nonumber
\end{eqnarray}
provided that $||f||_{L^1(I)} < \varepsilon$.
\end{corollaire}

\section{ Application}
In this section, we prove a logarithmic stability result for the inverse problem of identification of a Robin parameter in two dimensional annular domain.
Let $I$ be  any connected open subset of the boundary of the annular $G_s $ and let $ J = \partial G_s \setminus I$.
We consider the following inverse problem $(I. P)$.

Given a function $\varphi$ and a prescribed flux $\phi$  on $I$, find a function $q \in {\mathbf{Q}}_{ad}^n$ such that the solution $u$ to the problem
\begin{eqnarray}
(N. R)\left\{
\begin{array}{lrl}
 \ds \triangle u = 0 &\hbox {\rm in} &G_s,\\
 \ds\partial_n u = \Phi &\hbox{\rm on}
 &I,\\
  \ds\partial_n u + qu=0 &\hbox {\rm   on} & J,
\end{array}\right.\nonumber
\end{eqnarray}
also satisfies $u_{|I} = \varphi$,

where $\ds \partial_n $ stands for the partial derivative with respect to the outer normal unit vector to $\partial G_s$ and  the admissible set ${\mathbf{Q}}_{ad}^n$ of smooth Robin coefficient is defined by:
\begin{eqnarray}
 {\mathbf{Q}}_{ad}^n =\left\{ q \in {\mathcal{C}}^{n}_0({\overline J}),\ |q^{(k)}| \leq c',\ 0\leq k\leq n,\ \hbox{and}\ \ q \geq c  \right\},\nonumber
\end{eqnarray}

where $c, c' $ are non negative constants and $K$ is a nonempty connected subset of $J$  far from the boundary of $J$.
For $ q \in {\mathbf{Q}}_{ad}^n$, we denote by $ u_{q} $ the
solution of the  Neumann-Robin problem $ (N. R)$.

Referring to \cite{Chaabane Jaoua 99, chaabane jaoua leblond 2003,
chaabane ferchichi kunisch 2004}, we have the following

\begin{lemme}\label{identlemma2}(\cite{Chaabane Jaoua 99, chaabane jaoua leblond 2003,
chaabane ferchichi kunisch 2004})
 Let $n \in {\mathbb{N}}, \ \Phi \in W^{n,2}(I ) $ with non-negative value such that
 $\phi\ \not\equiv 0 $ and assume that $q \in {\mathbf{Q}}_{ad}^n$ for some constants $c, c' >0$.
 Then the solution $u_q$ of the inverse problem $(I. P)$ belongs to $W^{n+3/2,2}(G_s)$.

Furthermore, there exist non negative constants $ \alpha,\  \beta  $
such that for every $q \in {\mathbf{Q}}_{ad}^n$ and every $\Phi \in W^{n,2}(I )$, we have
\begin{eqnarray}
 u_q \geq \alpha > 0 \quad  \hbox{and}\quad  ||u||_{W^{n+1,2}(\partial G_s)} \leq \beta.\nonumber
\end{eqnarray}
\end{lemme}

The following identifiability result proves the uniqueness of the solution $q$ of the inverse problem $(I. P)$.

\begin{lemme}\label{identlemma} (\cite{Chaabane Jaoua 99, chaabane ferchichi kunisch
2004}) The mapping
\begin{eqnarray}
\begin{array}{lccl}
    F:&{\mathbf{Q}}^n_{ad}& \longrightarrow& L^2(\Gamma_d)\\
      &                q & \longmapsto & u_{q_{/_{\Gamma_d}}}
\end{array}\nonumber
\end{eqnarray}
is well defined, continuous and injective.
\end{lemme}

According to Theorem \ref{theorem:kinfty}, we establish the following
stability result.

\begin{theoreme}\label{theorem:stability}
Let $n \geq 2$ and  $\phi \in W_0^{n,2}(I)$ such that $\phi \not\equiv \, 0$ \
and \ $\phi \geq 0$. Then, there exists a non negative constant $C$
such that  for any $ q_1, q_2 \in {\mathbf{Q}}^n_{ad}$, we have
\begin{eqnarray}
 \|q_1 - q_2\|_{L^\infty(J)} \leq \frac{C}{\left|Log \|u_{q_1}
- u_{q_2}\|_{L^1(I)}\right|^{n-1}},\nonumber
\end{eqnarray}
provided that $\|u_{q_1} - u_{q_2}\|_{L^1(I)} < e^{-\frac{q_0}{\lambda C_s}}$.
\end{theoreme}
\medskip
{\bf Proof.} Referring to (\cite{Leblond}, Lemma 12), we introduce for every $i =1,2$, the analytic function $f_i$ in $G_s$ satisfying  $u_{q_i} = Re f_i$ and $ f_i \in H^{n+1,2}(\partial G_s)$. Moreover, Lemma  \ref{identlemma2} together with the Gagliardo-Nirenberg inequalities prove that there exists  non negative constants $M, K$ depending only on $s$ and the class ${\mathbf{Q}}^n_{ad}$ such that
\begin{eqnarray}\label{eqn:bounded}
||f_i||_{H^{n,\infty}(G_s)}\leq M ||f_i||_{H^{n+1,2}(G_s)} \leq K \quad \hbox{for}\ i=1,2
\end{eqnarray}
Using the equation $ \ds\partial_n u + q u = 0 \quad \hbox {\rm   on} \quad   J,$
we get for \ $f = f_1 -f_2$ \ that
\begin{eqnarray}
\ds q_1 - q_2 = - \frac{1}{Re f_1}\ \frac{\partial Im f_1}{\partial \theta} + \frac{1}{Re f_2}\ \frac{\partial Im f_2}{\partial \theta}\ =\ - \frac{1}{Re f_1}\ \frac{\partial Im f}{\partial \theta} + \frac{\partial Im f_2}{\partial \theta}\ \frac{Re f}{Re f_1 \ Re f_2}.\nonumber
\end{eqnarray}
It follows from Lemma \ref{identlemma2}, that
\begin{eqnarray}
 \|q_1 - q_2\|_{L^\infty(J)} \leq  \frac{1}{\alpha}\|f\|_{W^{1, \infty}(J)} +
\frac{\beta}{ \alpha^2} \  \|f\|_{L^\infty(J)}
\leq \ \left(\frac{1}{\alpha} + \frac{\beta}{ \alpha^2}\right) \|f\|_{W^{1, \infty}(J)}.\nonumber
\end{eqnarray}
Hence, from (\ref{eqn:bounded}) and Corollary \ref{cor:k-jinfty} we get
\begin{eqnarray}
 \|q_1 - q_2\|_{L^\infty(J)} \leq \frac{C}{\left|Log \|u_{q_1}
- u_{q_2}\|_{L^1(I)}\right|^{n-1}},\nonumber
\end{eqnarray}
provided that $\|u_{q_1} - u_{q_2}\|_{L^1(I)} < e^{-\frac{q_0}{\lambda C_s}}$.

The particular case where  $ I = {\TT}$, has been recently established by Leblond and al. in \cite{Leblond}.

\begin{corollaire}\label{cor:stability}
Let $n \geq 2$, let $\phi \in W_0^{n,2}({\TT})$\ such that $\phi \not\equiv \, 0$ \
and \ $\phi \geq 0$. Then, there exists a non negative constant $C$
such that  for any $ q_1, q_2 \in {\mathbf{Q}}^n_{ad}$, we have
\begin{eqnarray}
 \|q_1 - q_2\|_{L^\infty(s{\TT})} \leq \frac{C}{\left|Log \|u_{q_1}
- u_{q_2}\|_{L^1({\TT})}\right|^{n-1}},\nonumber
\end{eqnarray}
provided that $\|u_{q_1} - u_{q_2}\|_{L^1({\TT})} < e^{-\frac{q_0}{\lambda C_s}}$.
\end{corollaire}

{\small
{\em Authors' addresses}:
{\em Imed Feki}, Department of Mathematics, Faculty of Sciences,
Sfax University, B.P 1171, Sfax 3018, Tunisia.

 e-mail: \texttt{imed.feki@fss.rnu.tn}.

}

\end{document}